\documentclass[a4paper, 12pts]{article}

\usepackage{arxiv}

\usepackage[utf8]{inputenc} 
\usepackage[T1]{fontenc}    
\usepackage{hyperref}       
\usepackage{url}            
\usepackage{booktabs}       
\usepackage{amsfonts}       
\usepackage{nicefrac}       
\usepackage{microtype}      
\usepackage{lipsum}
\usepackage{graphicx}
\usepackage[utf8]{inputenc}
\usepackage{caption}
\usepackage{graphicx}
\usepackage{amsmath}
\usepackage{subcaption}
\usepackage[version=4]{mhchem}
\usepackage{siunitx}
\usepackage{booktabs}
\hypersetup{hidelinks}
\usepackage[graphicx]{realboxes}
\usepackage[table,xcdraw]{xcolor}
\usepackage[inkscapeformat=png]{svg}
\usepackage{float}
\usepackage{longtable,tabularx}
\usepackage{lineno}

\title{Data-driven Discovery of The Quadrotor Equations of Motion Via Sparse Identification of Nonlinear Dynamics}

\author{
 Zeyad M. Manaa\thanks{Corresponding author: Z.~M.~M \texttt{\{g202216800@kfupm.edu.sa\}}.\\
\textit{Author contributions}: Z.~M.~M formulated research; Z.~M.~M and M.~R.~E performed research and simulations; Z.~M.~M, M.~R.~E, and A.~M.~A analyzed data; M.~R.~E and Z.~M.~M wrote the manuscript; and Z.~M.~M, M.~R.~E, and A.~M.~A reviewed and finalized the manuscript.\\
Z.~M.~M would like to acknowledge the support provided by the Deanship of Research Oversight and Coordination at King Fahd University of Petroleum and Minerals (KFUPM) under Research Grant xxxxx.} \\
  Aerospace Engineering Department\\
  King Fahd University of Petroleum and Minerals\\
  Dhahran, 31261, Saudi Arabian\\
   \And
  Mohammed R. Elbalshy \\
  Aerospace Engineering Department\\
  King Fahd University of Petroleum and Minerals\\
  Dhahran, 31261, Saudi Arabian\\
  \And
 Ayman M. Abdallah \\
 Aerospace Engineering Department\\
 King Fahd University of Petroleum and Minerals\\
 Dhahran, 31261, Saudi Arabian\\
}

\begin{document}
\maketitle
\begin{abstract}
Dynamical systems provide a mathematical framework for understanding complex physical phenomena. The mathematical formulation of these systems plays a crucial role in numerous applications; however, it often proves to be quite intricate. Fortunately, data can be readily available through sensor measurements or numerical simulations. In this study, we employ the Sparse Identification of Nonlinear Dynamics (SINDy) algorithm to extract a mathematical model solely from data. The influence of the hyperparameter $\lambda$ on the sparsity of the identified dynamics is discussed. Additionally, we investigate the impact of data size and the time step between snapshots on the discovered model. To serve as a data source, a ground truth mathematical model was derived from the first principals, we focus on modeling the dynamics of a generic 6 Degrees of Freedom (DOF) quadrotor. For the scope of this initial manuscript and for simplicity and algorithm validation purposes, we specifically consider a sub-case of the 6 DOF system for simulation, restricting the quadrotor's motion to a 2-dimensional plane (i.e. 3 DOF). To evaluate the efficacy of the SINDy algorithm, we simulate three cases employing a Proportional-Derivative (PD) controller for the 3 DOF case including different trajectories. The performance of SINDy model is assessed through the evaluation of absolute error metrics and root mean squared error (RMSE). Interestingly, the predicted states exhibit at most a RMSE of order of magnitude approximately $10^{-4}$, manifestation of the algorithm's effectiveness. This research highlights the application of the SINDy algorithm in extracting the quadrotor mathematical models from data. We also try to investigate the effect of noisy measurements on the algorithm efficacy. The successful modeling of the 3 DOF quadrotor dynamics demonstrates the algorithm's potential, while the evaluation metrics validate its performance, thereby clearing the way for more applications in the realm of unmanned aerial vehicles.
\keywords{{Keywords: dynamical systems \and machine learning \and sparse regression \and system identification \and quadrotor \and numerical simulations \and optimization}}
\end{abstract}

\section{Introduction}
{Dynamical systems provide a mathematical representation for describing the world. Formally, dynamical systems are concerned with the analysis, and interpretation of the behavior of sets of differential equations that trace the advancement of a system's state across time \cite{arrowsmith_introduction_1990, brunton_data-driven_2019}. Classical dynamics has been discussed for years and applied to a wide range of applications; it is far too complex to be stated in a few words. However, the combination of big data, machine learning techniques, and statistical learning is driving a revolution in data-driven dynamical modeling and control of complex systems, with analytical derivations being replaced by data-driven methodology \cite{kaiser2018sparse}. In most real-world circumstances, data is so abundant that it cannot be comprehended. Furthermore, the physical principles that control these data are frequently complex; this is true for most physical concerns, such as climate science, epidemiology, and finance, to name a few. Researchers \cite{rowley2004model, rowley2017model, hayhoe2020data} attempt to decompose data in order to acquire insights into these massive, yet unexplained, datasets. As a result, the abundance of data remains a challenge because it is inherently imperfect. Consequently, the rise of data-driven dynamics is paving the way for such difficulties. 

The outcome of this development in data science was not initially noticeable on dynamical systems \cite{brunton_discovering_2016}, but recently a lot of research has been geared towards that area. Bongard and Lipson \cite{bongard2007automated} and Schmidt and Lipson \cite{schmidt2009distilling} was able to glean information about a nonlinear dynamical system's structure which was used after that to find the nonlinear differential equations \cite{koza1994genetic} by symbolic regression. Dynamic Mode Decomposition (DMD) was first introduced to fluid dynamics by Schmid \cite{schmid2008dynamic}. He extracted spatiotemporal structures from high dimensional data using DMD to analyse it. Researchers have extended the DMD method to explore other interesting intersections, including but not limited to the Koopman operator \cite{koopman_hamiltonian_1931}, extended DMD \cite{williams_datadriven_2015} to allow for approximation for Koopman operator, kernel DMD \cite{williams_kernel-based_2015, baddoo2022kernel}, time-delay DMD \cite{brunton_chaos_2017, arbabi_ergodic_2017}. Moreover, Proctor et al. \cite{proctor_dynamic_2016} extended the DMD to account for control input. These methods allow for linear system identification. Although these techniques rely on a precise set of coordinates for dynamics linearization, they are beneficial for developing dynamics-based linear models that progress high-dimensional observations across time  \cite{tu_dynamic_2014, tu2013dynamic} which are useful for control purposes \cite{champion2019data}. Recent studies \cite{yeung_learning_2017, takeishi_learning_2018} have examined the use of deep learning techniques to choose the proper set of coordinates to be used in both DMD and extended DMD.

While linear dynamical systems have several great qualities, there are many complex and intriguing dynamical phenomena that a linear model cannot properly capture. This encourages ways for discovering models of nonlinear dynamical systems. Consequently, another breakthrough on data-driven dynamics happened after Brunton et al. \cite{brunton_discovering_2016} revealed the method of discovering nonlinear equations of motion using SINDy algorithm. They combined compressed sensing \cite{candes_robust_2006, candes_stable_2006, donoho_compressed_2006} with sparse regression \cite{tibshirani_regression_1996, james_introduction_2021} and came up with the new technique of SINDy, which is based on the idea that most dynamical systems contain a small number of active terms that account for the majority of dynamics. Recently, SINDy algorithm has been generalized to include control inputs \cite{kaiser2018sparse}, to include tensor bases \cite{gelss2019multidimensional}, to discover partial differential equation \cite{rudy2017data, schaeffer2017learning}, and to account for noisy data as a kind of challenge \cite{schaeffer_sparse_2017, reinbold2020using}. The usage of SINDy algorithm thereafter begin to be a workhorse in the fields of optics, chemical engineering, robotics, and disease modeling \cite{jiang_modeling_2021, thiele_system_2020, bhadriraju_operable_2020, bhattacharya_sparse_2020, sorokina_sparse_2016}.

Despite the increasing body of literature in the field, there has been a dearth of research focused on the domain of aerospace, especially the Unmanned Aerial Vehicles (UAV). A few studies discussed data-driven dynamics and control as an application of the UAVs. Manzoor et al. \cite{manzoor2022model} proposed a novel approach for controlling ducted fan aerial vehicles (DFAVs) by integrating model predictive control (MPC) and physics-informed machine learning. They combined the physics-based model with the data-driven model to account for sudden changes in the system dynamics motivated by \cite{quade2018sparse, kaheman2019learning}. Also, Kaiser et al. \cite{kaiser2018sparse} proposed an integration between SINDy and MPC. They found that the SINDy method needs less data compared to other popular machine learning approaches, such as neural networks. This reduces the likelihood of over-fitting. In contrast to \cite{manzoor2022model}, when combined with MPC, the SINDy algorithm yields models that are both computationally feasible and precise, even when trained on limited data. 

Motivated by this, SINDy algorithm will be utilized to discover the EOM of quadrotor using numerical simulation data at first, and for the next phase, experimental data might be used. Moreover, a study on SINDy's sparsity hyperparameter effect on model discovery will be examined along with the effect of the sampling time of data collected whether numerically or experimentally. 

}

\section{Modelling and Simulation}
\subsection{Coordinate Systems}
\begin{figure}[t]
\centering
\includegraphics[width=.7\textwidth]{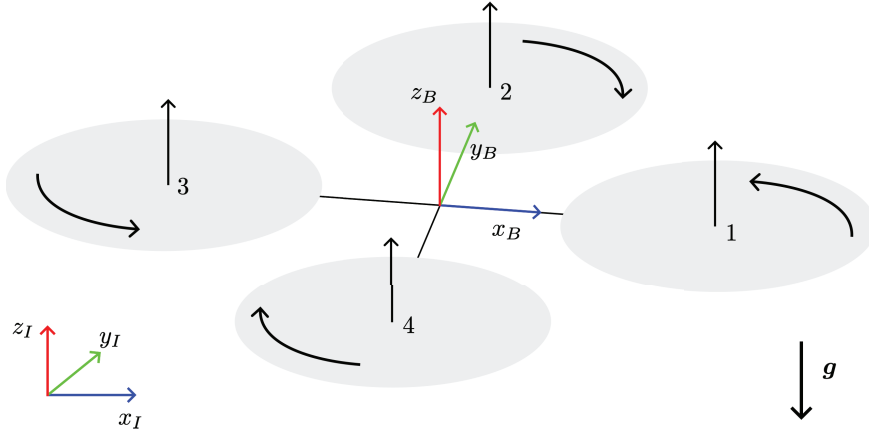}
\caption{Quadrotor with the body and the inertial reference frames.}
\label{fig:coordinate systems}
\end{figure}

The coordinate systems for the quadrotor are shown in Fig. \ref{fig:coordinate systems}. The inertial reference frame, $\mathcal{I}$, is defined by 
\begin{equation}
    \mathcal{I} = \{\mathcal{O}, \boldsymbol{x}_{\mathcal{I}}, \boldsymbol{y}_{\mathcal{I}}, \boldsymbol{z}_{\mathcal{I}}\}
\end{equation}
where, $\mathcal{O}$ is the frame's origin, $\mathbf{x}_{\mathcal{I}}$, $\mathbf{y}_{\mathcal{I}}$, $\mathbf{z}_{\mathcal{I}}$ are unit vectors defining the frame axes. 
The body frame, $\mathcal{B}$, is fixed to the center of mass of the quadrotor defined by
\begin{equation}
     \mathcal{B} = \{\mathcal{C}, \boldsymbol{x}_{\mathcal{B}}, \boldsymbol{y}_{\mathcal{B}}, \boldsymbol{z}_{\mathcal{B}}\}
\end{equation}
where, $\mathcal{C}$ is the frame's origin, $\mathbf{x}_{\mathcal{B}}$, $\mathbf{y}_{\mathcal{B}}$, $\mathbf{z}_{\mathcal{B}}$ are unit vectors defining the frame axes.

We may acquire the transformation matrix that transition between the body frame to the inertial frame by using the $\psi-\theta-\phi$ sequence, which indicates the yaw, pitch and roll angles respectively. 

\begin{equation}
{ }^{\mathcal{I}}\mathcal{R}_{\mathcal{B}}=\left[\begin{array}{ccc}
c \psi c \theta-s \phi s \psi s \theta & -c \phi s \psi & c \psi s \theta+c \theta s \phi s \psi \\
c \theta s \psi+c \psi s \phi s \theta & c \phi c \psi & s \psi s \theta-c \psi c \theta s \phi \\
-c \phi s \theta & s \phi & c \phi c \theta
\end{array}\right]
\end{equation}
Giving the fact that this matrix is orthonormal
\begin{equation}
{}^{\mathcal{I}}\mathcal{R}_{\mathcal{B}}{}^{\mathcal{I}}\mathcal{R}^{-1}_{\mathcal{B}}= {}^{\mathcal{I}}\mathcal{R}_{\mathcal{B}}{}^{\mathcal{I}}\mathcal{R}^{\top}_{\mathcal{B}} = \texttt{eye(3)}
\end{equation}
where \texttt{eye(3)} is a $3\times3$ identity matrix.

\subsection{Quadrotor Dynamics}
\label{sec:EOM}
We assume the quadrotor is a rigid body with 6 DOF, having a mass m and Intertia matrix $\boldsymbol{J} = \texttt{diag}(J_x, J_y, J_z)$.
Let $\boldsymbol{r}$ denotes the position of $\mathcal{C}$ in $\mathcal{I}$, the $\sum \boldsymbol{F}$ denotes the summation of forces acting upon the body and $T_n$ denotes the motor force. The equations governing the motion of $\mathcal{C}$ is derived by applying Newton's law to the translational motion considering the control actuation applied on the body coordinate.

\begin{equation}
    \label{eqn: acc_body}
\begin{aligned}   
    m\ddot{\boldsymbol{r}}^{\mathcal{I}}&= m\left(\ddot{\boldsymbol{r}}^{\mathcal{B}}+{\boldsymbol{\omega}^{\mathcal{B/I}}\times{\dot{\boldsymbol{r}}}^{\mathcal{B}}}\right) =\sum{\boldsymbol{F}}\\
     \sum{\boldsymbol{F}} &= {\boldsymbol{F}_g} + {\boldsymbol{F}}_{\text{th}} = \begin{bmatrix}
        0\\0\\-mg
    \end{bmatrix} + {}^{\mathcal{I}}\mathcal{R}_{\mathcal{B}} \begin{bmatrix}
        0\\0\\\sum_{n=1}^{4}{T_n}
    \end{bmatrix}\\
    \ddot{\boldsymbol{r}}^{\mathcal{B}} &= \frac{1}{m}\sum{\boldsymbol{F}} - {\boldsymbol{\omega}^{\mathcal{B/I}}\times{\dot{\boldsymbol{r}}}^{\mathcal{B}}}
\end{aligned}
\end{equation}

We also employ Euler’s equation to model the attitude dynamics such that $J$ is the quadrotor’s moment of inertia, $\boldsymbol{\omega}$ is the angular velocity, $\boldsymbol{h}$ is the angular momentum defined as $\boldsymbol{h} = \boldsymbol{J} \boldsymbol{\omega}$, the $\sum \boldsymbol{M}$ is the summation of the moments acting upon the system. 
\begin{equation}
\label{eqn: angular acc body}
    \dot{\boldsymbol{\omega}} = \boldsymbol{J}^{-1} \left({\sum \boldsymbol{M}} - \boldsymbol{\omega}\times{\boldsymbol{J}\boldsymbol{\omega}}\right)
\end{equation}
The relation between quadrotor velocity in the body frame and position in the inertial frame can be expressed as $\dot{\boldsymbol{r}}^{\mathcal{I}} = {}^{\mathcal{I}}\mathcal{R}_{\mathcal{B}} \dot{\boldsymbol{r}}^{\mathcal{B}}$. Also, to get the attitude of the quadrotor we define the relation between quadrotor compenets of angular velocity in the body frame and the roll, pitch and yaw derivatives in the inertial frame by $\boldsymbol{\omega}=\mathcal{T}\boldsymbol{\dot{a}}$. In summary, the overall dynamics of the system in hand are given by equation \ref{eqn: overall dynmaics}.

\begin{equation}
\label{eqn: overall dynmaics}
    \begin{aligned}
         \dot{\boldsymbol{r}}^{\mathcal{I}} & = {}^{\mathcal{I}}\mathcal{R}_{\mathcal{B}}\dot{\boldsymbol{r}}^{\mathcal{B}}\\
         \ddot{\boldsymbol{r}}^{\mathcal{B}} &= \frac{1}{m}\sum{\boldsymbol{F}} - {\boldsymbol{\omega}^{\mathcal{B/I}}\times{\dot{\boldsymbol{r}}}^{\mathcal{B}}}\\
        \boldsymbol{\omega} &=\mathcal{T}\boldsymbol{\dot{a}}\\
        \dot{\boldsymbol{\omega}} &= \boldsymbol{J}^{-1} \left({\boldsymbol{M}} - \boldsymbol{\omega}\times{\boldsymbol{J}\boldsymbol{\omega}}\right)
    \end{aligned}
\end{equation}

\subsection{Control inputs}
The first control actuation input is defined as the sum of the forces from each of the four rotors.
So,
\begin{equation}
\label{first input}
    u_1 = [T_1 + T_2 + T_3 + T_4]
\end{equation}
As shown in Fig \ref{fig:coordinate systems}, rotors one and three both rotate in the negative ${\boldsymbol{z}_{\mathcal{B}}}$ direction. Rotors two and four both rotate in the positive ${\boldsymbol{z}_{\mathcal{B}}}$ direction. Since the moment exerted on the quadrotor oppose the direction of blade rotation. So, the moments ${M_1}$ and ${M_2}$ are directed to the positive ${\boldsymbol{z}_{\mathcal{B}}}$ direction, while the moments ${M_3}$ and ${M_4}$ are directed to the positive ${\boldsymbol{z}_{\mathcal{B}}}$ direction.
Recalling to Eq. \ref{eqn: angular acc body}, it can be rewritten as
\begin{equation}
\label{eqn: second control input}
      {\boldsymbol{J}}\begin{bmatrix}
          \dot{p}\\
          \dot{q}\\
          \dot{r}\\
      \end{bmatrix} = \begin{bmatrix}
          {L(T_2-T_4)}\\
          {L(T_3-T_1)}\\
          {M_1-M_2+M_3-M_4}
      \end{bmatrix} - \begin{bmatrix}
          p\\
          q\\
          r\\
      \end{bmatrix}\times{\boldsymbol{J}}\begin{bmatrix}
          p\\
          q\\
          r\\
      \end{bmatrix}
\end{equation}
From Eq. \ref{first input} and Eq. \ref{eqn: second control input}, the total control actuation inputs, defined as the vecotr $\boldsymbol{U}=[u_1, u_2, u_3, u_4]^\top$, 
\begin{equation}
    \boldsymbol{U} = \begin{bmatrix}
        u_1\\
        u_2\\
        u_3\\
        u_4\\
    \end{bmatrix} = \begin{bmatrix}
          1 & 1 & 1 & 1\\
          0 & L & 0 & -L\\
          -L & 0 & L & 0\\
          \gamma & -\gamma & \gamma & -\gamma \end{bmatrix} \begin{bmatrix}
              T_1\\
              T_2\\
              T_3\\
              T_4\\
      \end{bmatrix}
\end{equation}
where ${L}$ is the distance from the rotor axis of rotation to the quadrotor center of mass and $\gamma=\frac{K_M}{K_F}$, where ${K_m}$ and ${K_t}$ are respectively the aerodynamic motor moment and force constants \cite[Section~3.2]{elkholy2014dynamic}.
\subsection{Simulation}
\begin{figure}
    \centering
    \includegraphics[width = 0.7\textwidth]{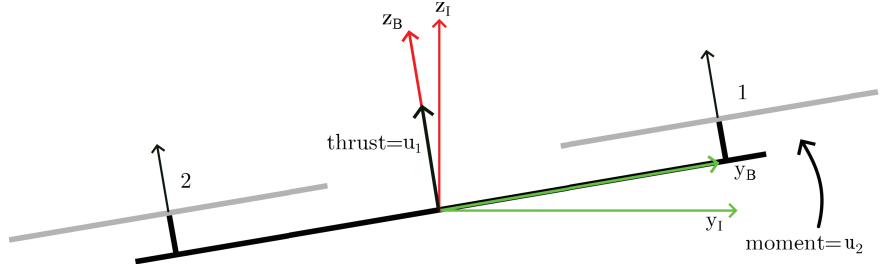}
    \caption{3 DOF quadrotor with the inertial and the body reference frames.}
    \label{fig:3DOF Quad}
\end{figure}

For the scope of this manuscript, we will partially focus on the simulation of a 3 DOF quadrotor system, as shown in Figure \ref{fig:3DOF Quad} Quad. This design, also known as a planar quadrotor, restricts the quadrotor's motion to the $y-z$ plane exclusively. At first, we wanted to use a simplified model to validate and test SINDy algorithm to the problem in hand. Then, gradually we transition to the entire 6-DOF analysis by studying this simplified subset derived from the comprehensive 6-DOF scenario mentioned in Section \ref{sec:EOM}. The simulation of the system's equations of motion will be done discretely using the well-known Runge-Kutta fourth-order (RK4) numerical technique.
\begin{equation}
    \boldsymbol{x}_{k+1} = \boldsymbol{f}_{\text{RK4}}(\boldsymbol{x}_k, \boldsymbol{u}_k, \Delta t)
\end{equation}
where $k$ is discrete time index, $\Delta t$ is the time step, $\boldsymbol{u}$ is the control actuation. For the full picture formulation of RK4 readers may refer to \cite{sola2017quaternion}. 

 The 3 DOF equations of motion are 
\begin{equation}
    \begin{aligned}
        \ddot{y} &= -\frac{u_1}{m} \sin{\phi} \\
        \ddot{z} &= \frac{u_1}{m} \sin{\phi} - g \\
        \ddot{\phi} &= \frac{u_2}{J_x}
    \end{aligned}
\end{equation}

or in a matrix format as
\begin{equation}
\label{eqn:3DOF Quad}
\left[\begin{array}{c}
\ddot{y} \\
\ddot{z} \\
\ddot{\phi}
\end{array}\right]=\left[\begin{array}{c}
0 \\
-g \\
0
\end{array}\right]+\left[\begin{array}{cc}
-\frac{1}{m} \sin (\phi) & 0 \\
\frac{1}{m} \cos (\phi) & 0 \\
0 & \frac{1}{I_{x x}}
\end{array}\right]\left[\begin{array}{l}
u_1 \\
u_2
\end{array}\right]
\end{equation}
We consider the state vector $\boldsymbol{x} = \begin{bmatrix}
    y, z, \phi, \dot{y}, \dot{z}, \dot{\phi}
\end{bmatrix}^{\top}$. Then, Equations \ref{eqn:3DOF Quad} are linearized about the hovering position using the first order Taylor series expansion approximation of the non-linear terms and developed PD controller for the robot.

The linearized version of Eqs. \ref{eqn:3DOF Quad} 
\begin{equation}
\begin{aligned}
\ddot{y} & =-g \phi \\
\ddot{z} & =-g+\frac{u_1}{m} \\
\ddot{\phi} & =\frac{u_2}{I_{x x}}
\end{aligned}
\end{equation}

Introducing a generic state variable $\boldsymbol{\varrho}$ that can represent either $\boldsymbol{y}$, $\boldsymbol{z}$, or $\boldsymbol{\phi}$ at a time. For this generic state variable $\boldsymbol{\varrho}$  to be driven to a new desired state, the commanded acceleration $\boldsymbol{\ddot{\varrho}}_c$ is needed. For this we define the proportional and derivative error as 
\begin{equation}
\begin{aligned}
& e_p=\varrho_d -\varrho \\
& e_v=\dot{\varrho}_d-\dot{\varrho}
\end{aligned}
\end{equation}
where the $\varrho_d$ is the desired value. Then we want both $e_p$ and $e_v$ to satisfy,
\begin{equation}
    \left(\ddot{\varrho}_d -\ddot{\varrho}_c\right)+k_p e_p+k_v e_v=0
\end{equation}
in order to guarantee the error's convergence under some values of $k_p$ and $k_v$.
So
\begin{equation}
\begin{aligned}
& u_1=m g+m \ddot{z}_c=m\left(\ddot{z}_d +k_{v, z}\left(\dot{z}_d -\dot{z}\right)+k_{p, z}\left(z_d-z\right) + g\right) \\
& u_2=J_{x} \ddot{\phi}_d=J_{ x}\left(\ddot{\phi}_c+k_{v, \phi}\left(\dot{\phi}_c-\dot{\phi}\right)+k_{p, \phi}\left(\phi_c-\phi\right)\right) \\
& \phi_c=-\frac{\ddot{y}_c}{g}=-\frac{1}{g}\left(\ddot{y}_d+k_{v, y}\left(\dot{y}_d-\dot{y}\right)+k_{p, y}\left(y_d-y\right)\right)
\end{aligned}
\end{equation}
From this, we can simulate the 3 DOF quadrotor based on different desired trajectory cases found in Table \ref{tab:cases}. 

We have set the \textit{mass} of the quadrotor to be $0.18$ Kg, the \textit{arm length} is 0.086 m, and the quadrotor's \textit{moment of inertia} ($J_x$) is 0.00025 Kg.m$^2$.

\begin{figure}
    \centering
    \includegraphics[width = \textwidth]{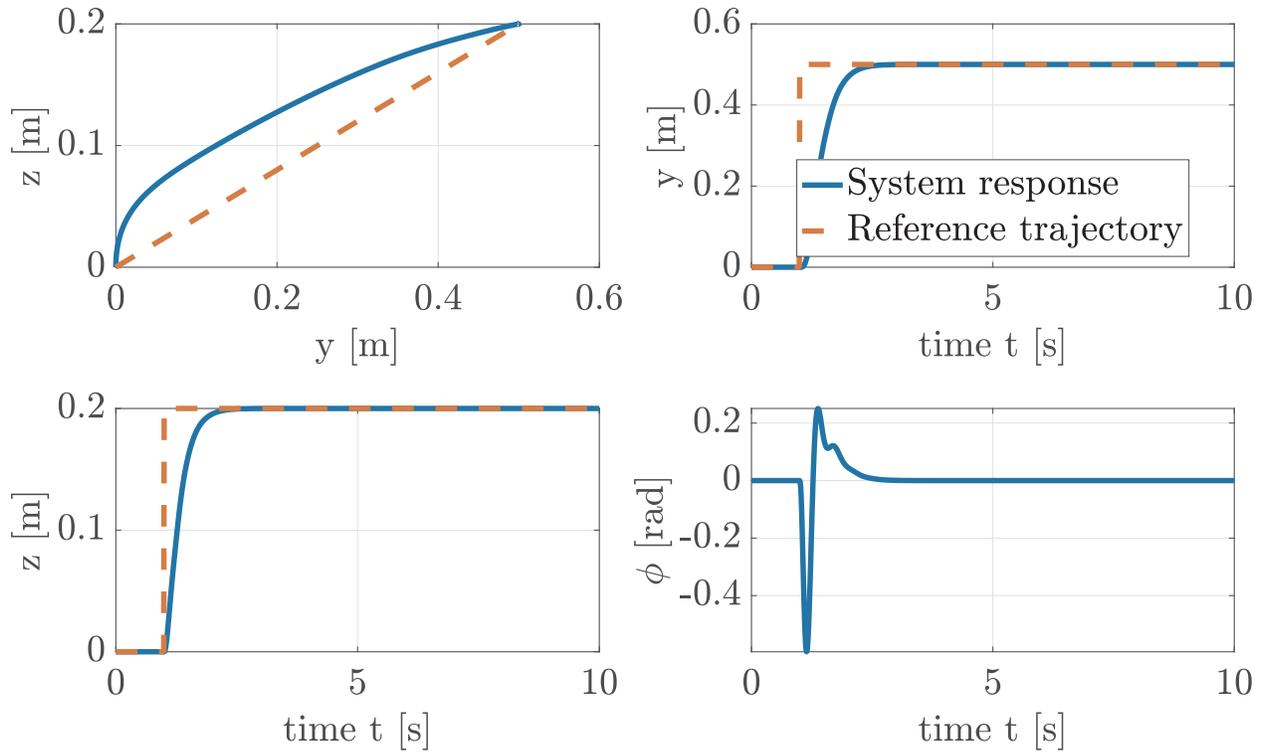}
    \caption{The reference step trajectory represent case A and the system response}
    \label{fig:response1}
\end{figure}

\begin{figure}
    \centering
    \includegraphics[width = \textwidth]{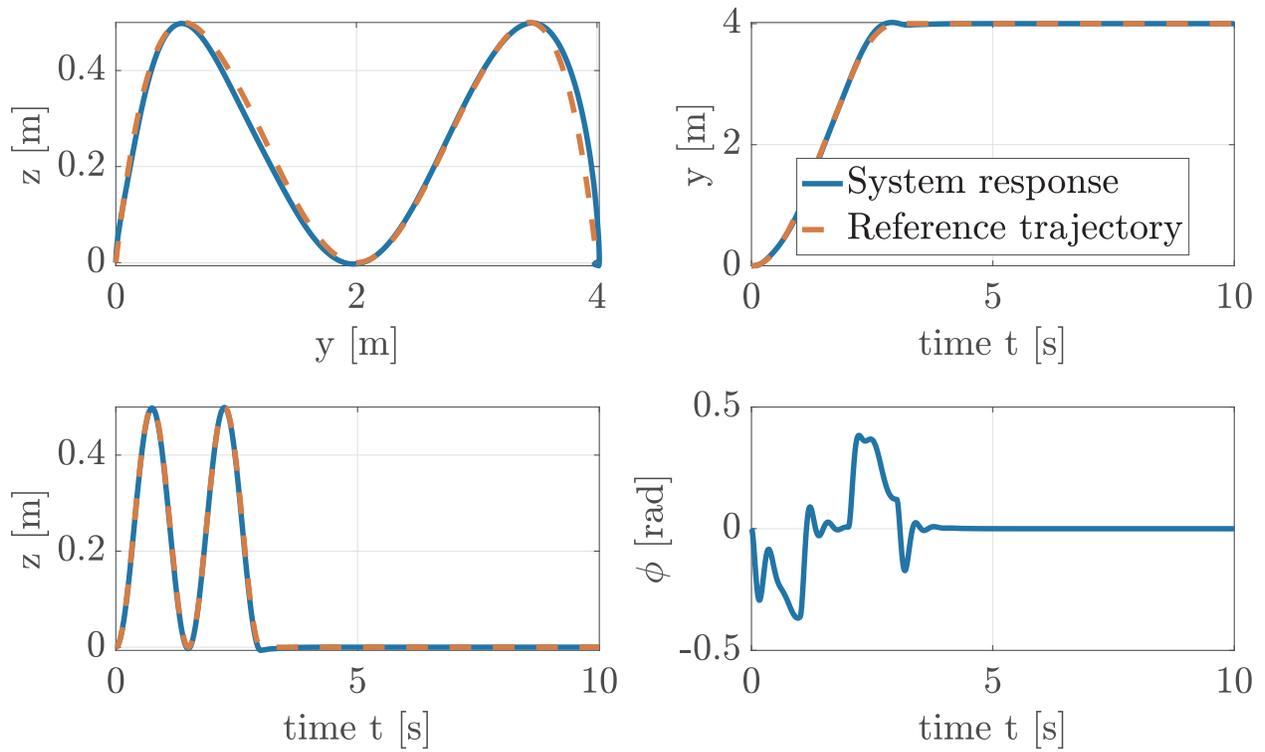}
    \caption{The reference sinusoidal trajectory represent case B and the system response}
    \label{fig:response2}
\end{figure}

\begin{figure}
    \centering
    \includegraphics[width = \textwidth]{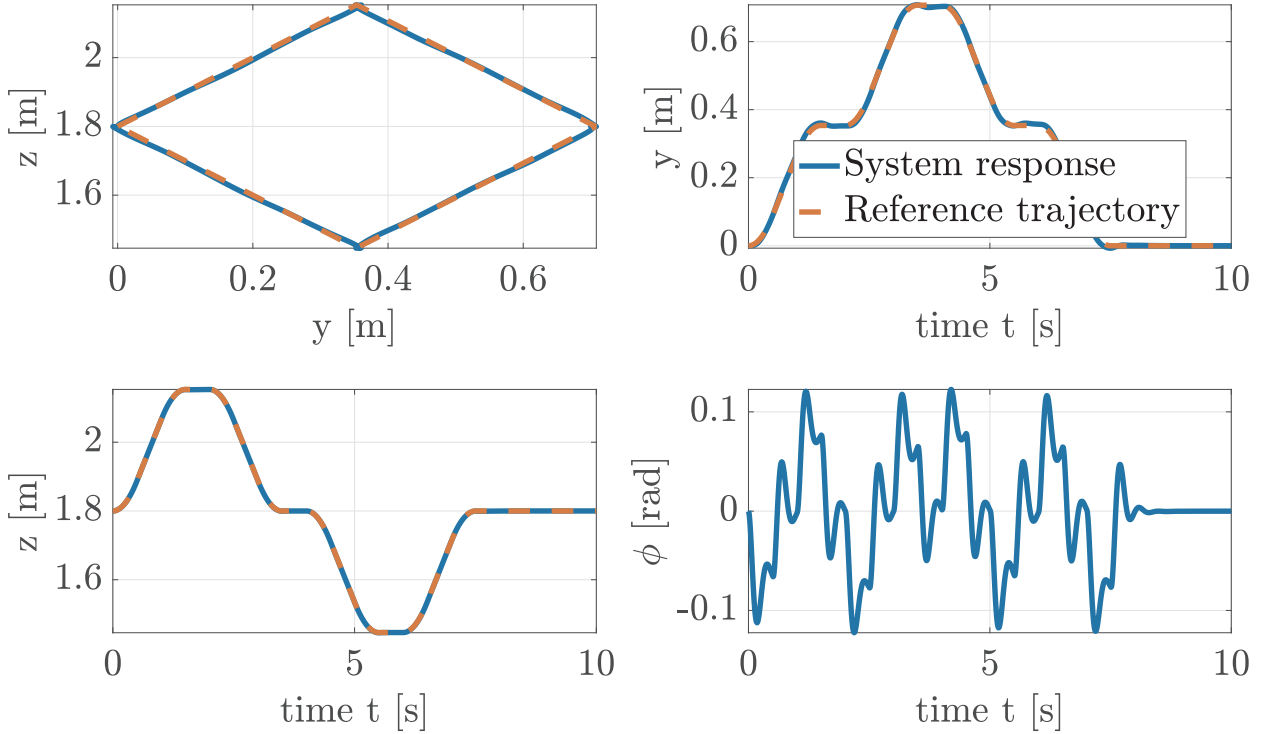}
    \caption{The reference diamond-shaped trajectory represent case C and the system response}
    \label{fig:response3}
\end{figure}

\begin{table}[t]
\centering
\begin{tabular}{@{}llll@{}}
\toprule
\textbf{Case} & \textbf{Trajectory type} & \textbf{Initial states} & \textbf{Desired states}\\ \midrule
A    & step            &   $x_0 = [0,\; 0,\;  0,\;  0,\; 0,\; 0]^{\top}$ & $x_d = [0.5,\; 0.2,\;  0,\;  0,\; 0,\; 0]^{\top}$               \\
B    & Sine            &   $x_0 =[0,\; 0,\;         0,\;  0,\; 0,\; 0]^{\top}$ & $x_d = [4,\; 0,\;  0,\;  0,\; 0,\; 0]^{\top}$                   \\
C    & Diamond         &   $x_0 =[0,\; 1.8,\;         0,\;  0,\; 0,\; 0]^{\top}$ & $x_d =[0,\; 0,\;         0,\;  0,\; 0,\; 0]^{\top}$                 \\ \bottomrule
\end{tabular}
\caption{Simulation cases for different trajectories.}
\label{tab:cases}
\end{table}
The response of the three cases is summarized in Figures \ref{fig:response1}, \ref{fig:response2}, and \ref{fig:response3}. In the simulation settings, one can see the acceptable tracking behaviour of the quadrotor. The resulted data from will be used as a training and test data for SINDy algorithm; more on this in section \ref{methodology}.

\section{Methodology}
\subsection{SINDy With Control Framework} \label{methodology}
Nonlinear dynamical systems can be represented as
\begin{equation} 
    \dot{\boldsymbol{x}} = \boldsymbol{f}(\boldsymbol{x})
\end{equation}
and we consider more general case when control takes place in the nonlinear dynamics described as
\begin{equation}\label{eqn:2}
    \dot{\boldsymbol{x}} = \boldsymbol{f}(\boldsymbol{x}, \boldsymbol{u})
\end{equation}
where $\boldsymbol{x} \in \mathbb{R}^n$ is the system states vector, $\boldsymbol{u} \in \mathbb{R}^q$ is the control vector, and $\boldsymbol{f}$ is a nonlinear such that $\boldsymbol{f}: \mathbb{R}^n \times \mathbb{R}^q \longrightarrow \mathbb{R}^n$.

SINDy with control technique utilizes sparse regression to identify a minimal set of active terms from a candidate library $\boldsymbol{\Pi(\boldsymbol{x}, \boldsymbol{u})}$ of linear and nonlinear terms in the system states $\boldsymbol{x}$ and actuation variables $\boldsymbol{u}$ that approximate the underlying nonlinear dynamics represented by $\boldsymbol{f}$. On the premise that a considerable number of systems manifest relatively sparse active terms in their dynamics. To construct $\boldsymbol{\Pi}$ we measure $m$ snapshots of the state $\boldsymbol{x}$ and actuation $\boldsymbol{u}$ by experiment or from numerical simulation. 
\begin{equation}
    \boldsymbol{X} = \begin{bmatrix}
\boldsymbol{x}^{\top}(t_1) \\
\boldsymbol{x}^{\top}(t_2) \\
\boldsymbol{x}^{\top}(t_3) \\
\vdots \\
\boldsymbol{x}^{\top}(t_m)
\end{bmatrix} = \begin{bmatrix}
x_1(t_1) & x_2(t_1) & \dots & x_n(t_1) \\
x_1(t_2) & x_2(t_2) & \dots & x_n(t_2) \\
x_1(t_3) & x_2(t_3) & \dots & x_n(t_3) \\
\vdots & \vdots & \ddots & \vdots \\
x_1(t_m) & x_2(t_m) & \dots & x_n(t_m)
\end{bmatrix}, \quad \boldsymbol{U} = \begin{bmatrix}
\boldsymbol{u}^{\top}(t_1) \\
\boldsymbol{u}^{\top}(t_2) \\
\boldsymbol{u}^{\top}(t_3) \\
\vdots \\
\boldsymbol{u}^{\top}(t_m)
\end{bmatrix} = \begin{bmatrix}
u_1(t_1) & u_2(t_1) & \dots & u_q(t_1) \\
u_1(t_2) & u_2(t_2) & \dots & u_q(t_2) \\
u_1(t_3) & u_2(t_3) & \dots & u_q(t_3) \\
\vdots & \vdots & \ddots & \vdots \\
u_1(t_m) & u_2(t_m) & \dots & u_q(t_m)
\end{bmatrix}
\end{equation}
The matrix $\boldsymbol{\Pi}$ can be reconstructed now by 
\begin{equation}
    \boldsymbol{\Pi}(\boldsymbol{X, U}) = \begin{bmatrix}
| & | & | & | & | & | &   & | & | &   \\
\boldsymbol{1} & \boldsymbol{X} & \boldsymbol{U} & \boldsymbol{X\otimes X} & \boldsymbol{X\otimes U} & \boldsymbol{U\otimes U} & \dots & \cos(\boldsymbol{X}) & \cos(\boldsymbol{X\otimes X}) & \dots \\
| & | & | & | & | & | &   & | & | &  
\end{bmatrix}
\end{equation}
The effectiveness of the candidate term library is important to the SINDy with control algorithm. A basic strategy starts with a simple option, like polynomials, and gradually increases the complexity of the library by incorporating additional terms \cite{kaiser2018sparse}. After evaluating the library, the states derivatives also should be evaluated. Here we evaluated the derivatives by using numerical approximation. Specifically, we used the finite difference method.

The system in Eq. (\ref{eqn:2}) can be rewritten as
\begin{equation}
    \dot{\boldsymbol{X}} = \boldsymbol{\Pi}(\boldsymbol{X}, \boldsymbol{U})\boldsymbol{\Omega}
\end{equation}
where $\boldsymbol{\Lambda}$ is a coefficient matrix that is almost sparse. This matrix tries to activate the fewest active terms in the candidate matrix $\boldsymbol{\Pi}$ that results in the best model fit:
\begin{equation}
    \omega_{j} = \operatorname*{argmin}_{\Tilde{\omega}_{j}} ||\dot{\boldsymbol{X}} - \boldsymbol{\Pi}(\boldsymbol{X, U})\Tilde{\omega}_{j}||_{2}^{2} + \lambda ||\Tilde{\omega}_{j}||_{1}
\end{equation}
This problem can be solved using variety of regression techniques. Including but not limited to, LASSO \cite{tibshirani_regression_1996}, sequential thresholded least-squares \cite{brunton_discovering_2016}, and Sparse Relaxed Regularized Regression (SR3) \cite{zheng2018unified}. We tried different optimizers and have found SR3 is the best for our case. 

\subsubsection{Hyperparameter Tuning}
The hyperparameter $\lambda$ is critical in identifying the most sparse dynamics. Typically, $\lambda$ takes values between 0 to 1. So, we scanned over the lambda domain with step size of 0.05 and evaluated the model using test data to cross validate utilizing the root mean squared error (RMSE). We found the best value of $\lambda$ to be equal to 0.45.

We assume the generic prediction of a generic state vector $\boldsymbol{\varrho}$ is denoted as $\hat{\boldsymbol{\varrho}}$ and the RMSE of $\boldsymbol{\varrho}$ is denoted as $\Tilde{\varrho}$
\begin{equation}
    \Tilde{\varrho} = \sqrt{\frac{\sum_{i=1}^{N} (\varrho_i - \hat{\varrho}_i)^2}{N}}
\end{equation}

The chosen lambda corresponds to the RMSE of the states of the planar quad case.

\begin{table}[htp!]
\centering
\label{tab:RMSE}
\begin{tabular}{@{}lllllll@{}}
\toprule
State                                                 & ${y}$ & ${z}$ & ${\phi}$ & $\dot{{y}}$ & $\dot{{z}}$ & $\dot{{\phi}}$ \\ \midrule
RMSE $\times 10^{-3}$ & 0.0088      & 0.0152      & 0.0227         & 0.0294            & 0.0070            & 0.1379               \\ \bottomrule
\end{tabular}
\label{RMSE}
\caption{The RMSE for the chosen $\lambda = 0.45$ for diamond shaped trajectory.}
\end{table}

\subsubsection{Candidate functions}
In the present study, we exploit our comprehensive understanding of the physical system in hand. Specifically, we propose that the presence of nonlinearity in the system can be expressed through polynomials and Fourier basis functions, such as those involving sine and cosine. Through mathematical analysis of the system, we have discerned that from the entire state space of the system, only the Euler angles can be represented as Fourier basis functions, and the rest can be characterized by polynomials.

\subsection{Training}
\label{sec:training}
The numerical data from the 3 DOF planar quadrotor simulation was utilized. We first choose Case C because we thought it would allow the quadrotor to span the entire state space. Consequently, it will be the ideal case to start with for training. It was subsequently demonstrated by the other cases that the algorithm is generic regardless of trajectory. As long as the trajectory gave sufficient data and spanned the state space appropriately, the algorithm successfully captured the complete dynamics. If, on the other hand, the quadrotor followed a signal that led it to move only in a straight line along the $y$ or $z$ axis, the algorithm may fail to distinguish the unseen states, resulting in an incomplete identification of the dynamics. We used snapshots with $\Delta t = 0.05$ with 1000 snapshots in time. We differentiated the states using finite difference numerical scheme of order one.

\section{Results and Discussion}
\subsection{Discovered Model}
The model is trained with the data extracted from case C as discussed earlier and it came out with the following dynamics
\begin{equation}
\begin{aligned}
    \dot{y} &= \dot{y} \\
    \dot{z} &= \dot{z} \\
    \dot{\phi} &= 0.993 \dot{\phi}\\
    \ddot{y} &= -5.549 u_1 \sin(\phi) \\
    \ddot{z} &= -9.811 + 5.556 u_1 \cos(\phi) \\
    \ddot{\phi} &= 4000.000 u_2 
\end{aligned}
\end{equation}
That nearly matches the original derived mathematical model from the first principles. Table \ref{tab:comp} shows a comparison between the discovered dynamics by SINDy and the gound truth mathematical model.
\begin{table}[]
\centering
\caption{Comparison between the discovered dynamics and ground truth mathematical model.}
\label{tab:comp}
\begin{tabular}{@{}lll@{}}
\toprule
\textbf{States} & \textbf{SINDy} & \textbf{Mathematical Model} \\ \midrule
$\dot{y}$ & {\color[HTML]{FE0000} $1.000 \dot{y}$} & {\color[HTML]{3166FF} $1.000 \;\dot{y}$} \\
$\dot{z}$ & {\color[HTML]{FE0000} $1.000\; \dot{z}$} & {\color[HTML]{3166FF} $1.000\; \dot{z}$} \\
$\dot{\phi}$ & {\color[HTML]{FE0000} $0.993 \;\dot{\phi}$} & {\color[HTML]{3166FF} 1.000 $\dot{\phi}$} \\
$\ddot{y}$ & {\color[HTML]{FE0000} $-5.549 \;u_1 \sin(\phi)$} & {\color[HTML]{3166FF} $- 5.5556\; u_1 \;\sin(\phi)$} \\
$\ddot{z}$ & {\color[HTML]{FE0000} $-9.811 + 5.556\; u_1 \cos(\phi)$} & {\color[HTML]{3166FF} $-9.81 + 5.556\; u_1 \cos(\phi)$} \\
$\ddot{\phi}$ & {\color[HTML]{FE0000} $4000.000 \; u_2$} & {\color[HTML]{3166FF} $4000.000 \; u_2$} \\ \bottomrule
\end{tabular}
\end{table}
However, we also used the other cases to train the model and it resulted in models that compare to the original one as we discussed in section \ref{sec:training}. 

\subsection{Testing}
\subsubsection{Case A}
\begin{figure}
     \centering
     \begin{subfigure}[b]{\textwidth}
         \centering
         \includegraphics[width=\textwidth]{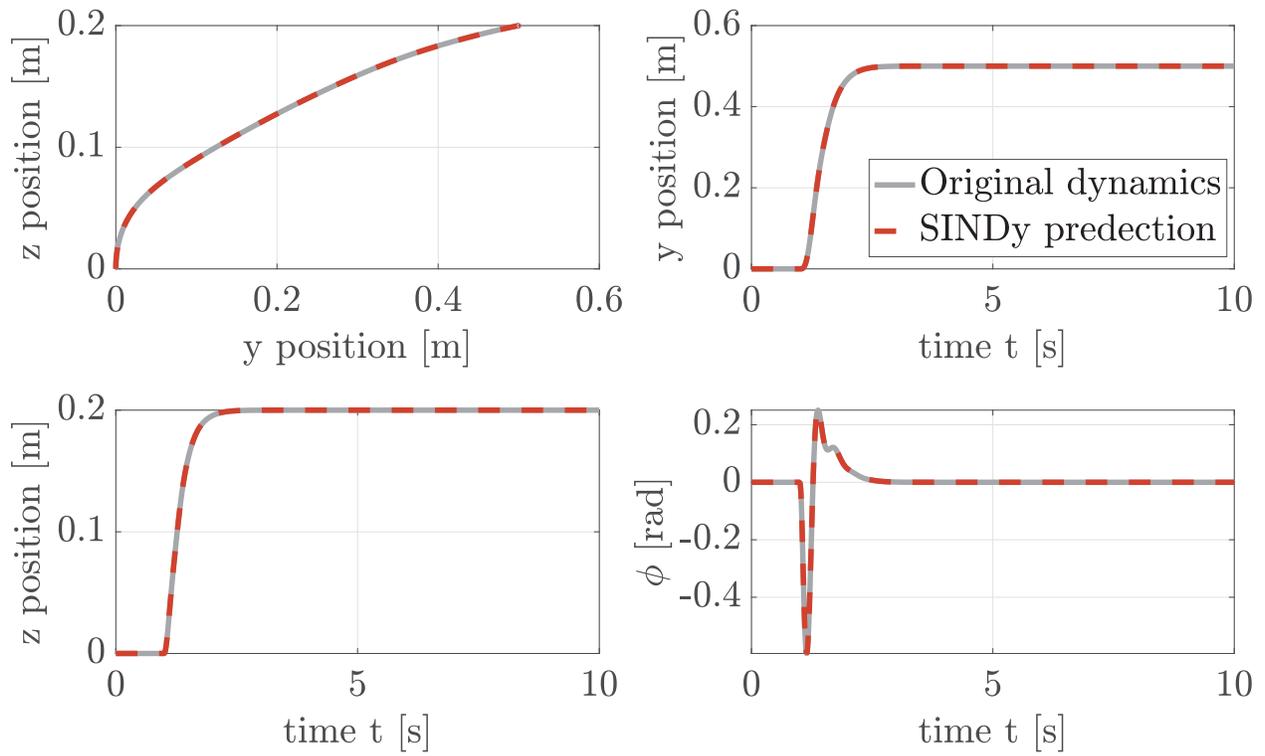}
         \caption{The reference step trajectory represent case A and the system response}
         \label{fig:y equals x}
     \end{subfigure}
     \hfill
     \begin{subfigure}[b]{\textwidth}
         \centering
         \includegraphics[width=\textwidth]{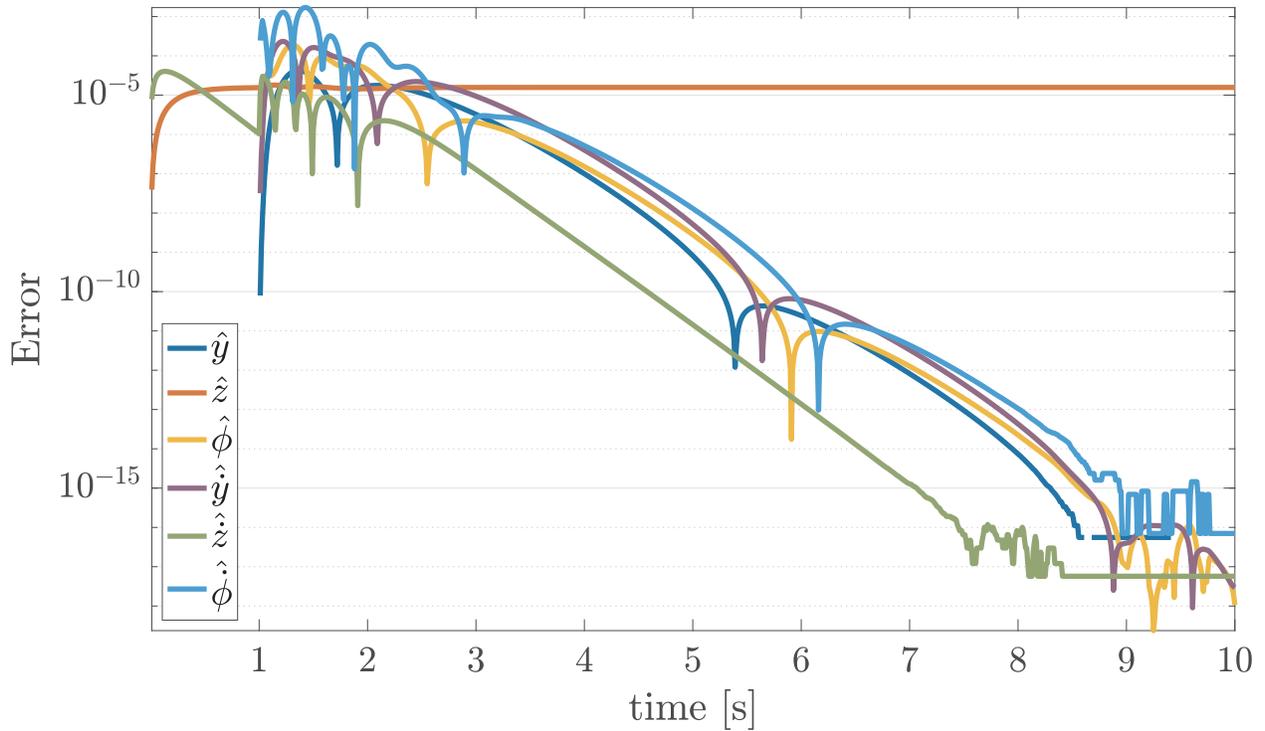}
         \caption{The absolute error between the predicted states $\hat{\boldsymbol{\varrho}}$ and the original states $\boldsymbol{\varrho}$}
         \label{fig:error_ramp}
     \end{subfigure}
     \label{fig:ramp_comp}
     \caption{SINDy comparison with the mathematical model.}
\end{figure}

Here we simulate the discovered dynamics using the step trajectory as desired trajectory. We compare the system behaviour over time for both the discovered SINDy dynamics and the ground truth. Figure \ref{fig:error_ramp} shows the absolute error between the predicted states $\hat{\boldsymbol{\varrho}}$ and the original states $\boldsymbol{\varrho}$. The results show that the absolute error approaches zero, giving strong validation for the identified model. This shows that the SINDy algorithm accurately captures the underlying dynamics of the system, closely matching the ground truth dynamics. The relatively small error between anticipated and original states confirms the discovered model's efficacy and reliability in capturing the main characteristics of quadrotor dynamics. The close agreement between the identified dynamics and the ground truth confirms the SINDy algorithm's utility as a good tool for extracting mathematical models from data.

\subsubsection{Case B}
\begin{figure}
     \centering
     \begin{subfigure}[b]{\textwidth}
         \centering
         \includegraphics[width=\textwidth]{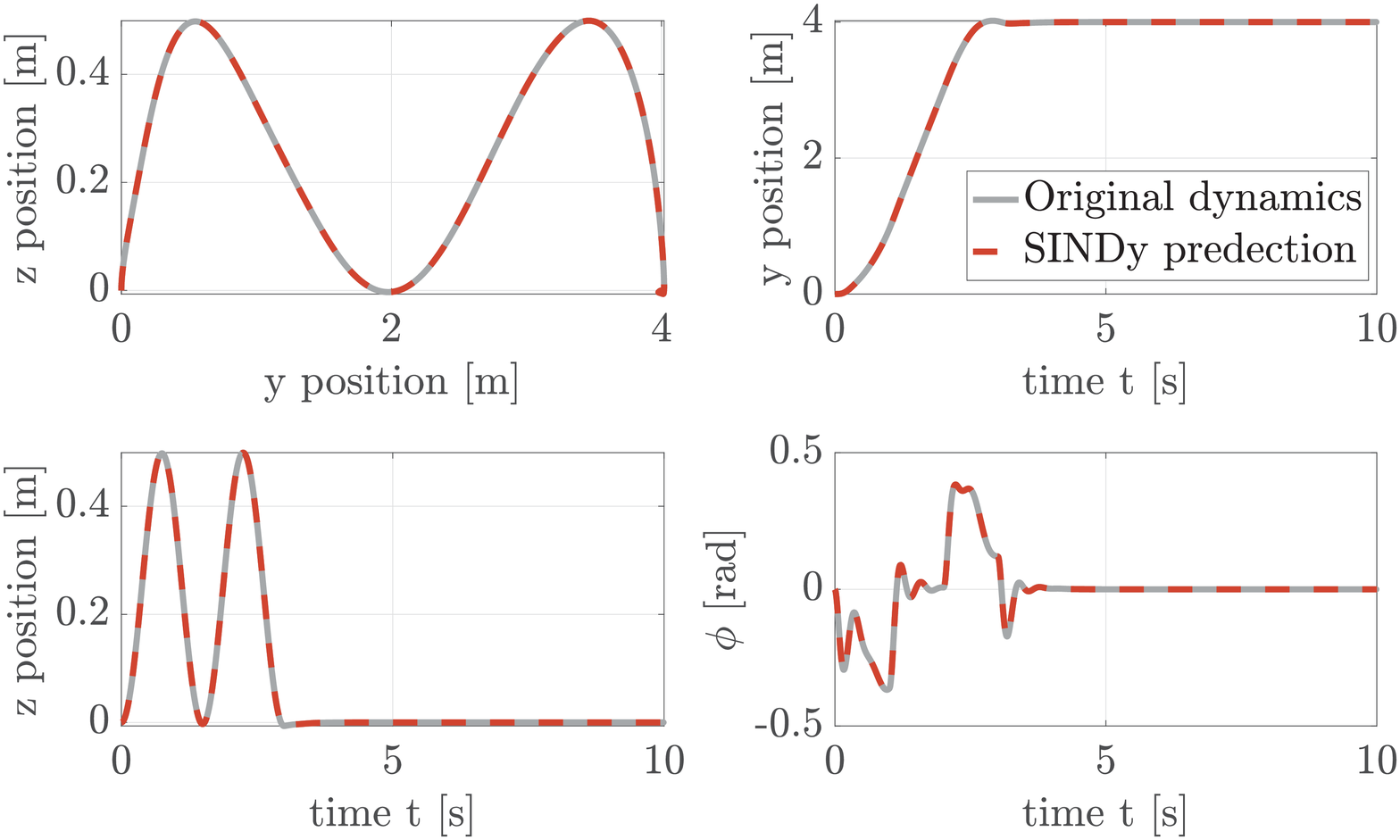}
         \caption{The reference sinusoidal trajectory represent case A and the system response}
         \label{fig:y equals x}
     \end{subfigure}
     \hfill
     \begin{subfigure}[b]{\textwidth}
         \centering
         \includegraphics[width=\textwidth]{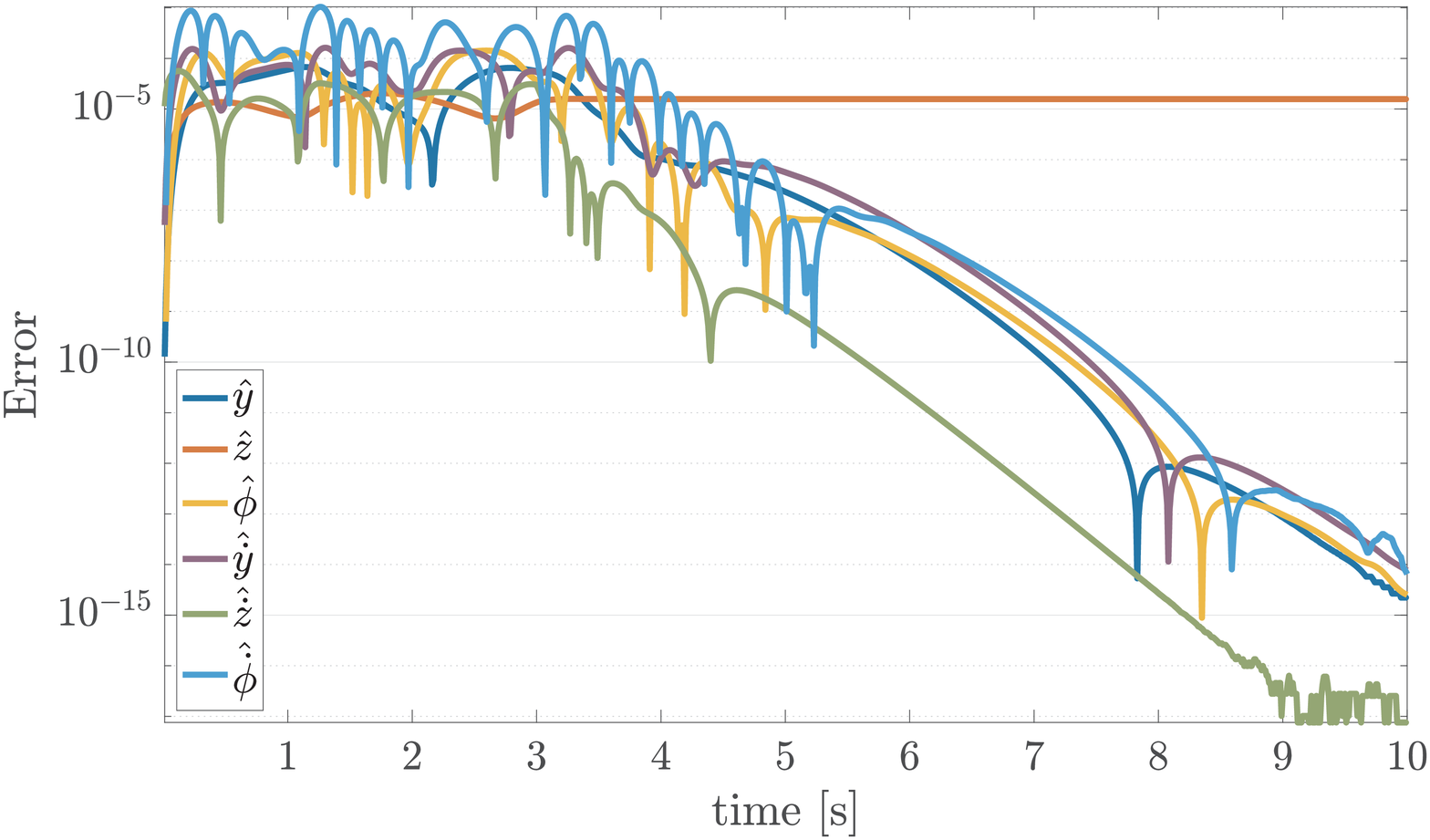}
         \caption{The absolute error between the predicted states $\hat{\boldsymbol{\varrho}}$ and the original states $\boldsymbol{\varrho}$}
         \label{fig:error_sin}
     \end{subfigure}
     \caption{SINDy comparison with the mathematical model.}
\end{figure}

As in the previous section, we simulate the discovered dynamics using the sinusoidal trajectory as desired trajectory. We compare the system behaviour over time for both the discovered SINDy dynamics and the ground truth. Figure \ref{fig:error_sin} shows the absolute error between the predicted states $\hat{\boldsymbol{\varrho}}$ and the original states $\boldsymbol{\varrho}$. We can say that the error nearly equal zero. This validates the discovered model. 

\section{Conclusion}
In this study, we have employed the SINDy algorithm to extract a mathematical model from simulation data of a quadrotor. Initially, the quadrotor dynamics were modeled in the 6 DOF configuration. Subsequently, for the purpose of simplification and model validation, the system was constrained to a 3 DOF representation within a 2-dimensional plane. To assess the effectiveness of the SINDy algorithm, three distinct simulation cases, as depicted in Figure \ref{fig:response}, were considered. A systematic exploration of the hyperparameter $\lambda$ space was conducted. Through this thorough analysis, we successfully identified the optimal model with an associated $\lambda$ value of 0.45, as demonstrated in Table \ref{tab:RMSE}, based on the RMSE metric. Table \ref{tab:comp} and Figures \ref{fig:error_ramp}, \ref{fig:error_sin} show that the algorithm captured almost the same dynamics as the quadrotor ground truth mathematical model.

We also demonstrated the SINDy powerfulness in obtaining solid mathematical models from data, which will aid in the advancement of modeling in the field of aerospace applications, particularly UAVs. Furthermore, in a future version of this research, we will attempt to analyze the algorithm's resilience in the face of noisy measurements in order to assure its robustness in real-world circumstances. Understanding how measurement errors affect the identified model's accuracy and reliability. Given the occurrence of mathematically unmodeled disturbances and uncertainties in real-world scenarios, future research should concentrate on building discrepancy models capable of properly capturing and accounting for these aspects. Incorporating these models with the SINDy algorithm can improve its ability to deal with disturbances, resulting in more robust and accurate predictions.

Furthermore, we intend to extend the SINDy algorithm's applicability to the full 6 DOF quadrotor model, encompassing all elements of its dynamic behavior. A more thorough understanding of the quadrotor's dynamics can be obtained by including the entire model, including translational and rotational motion. There will also be an in-depth examination of the optimizers used in the SINDy algorithm. Examining various optimization strategies and their impact on the quality of the identified model can result in improvements in convergence speed, accuracy, and the capacity to handle complicated datasets successfully.


\bibliographystyle{unsrt}  
\bibliography{references}


\end{document}